\documentclass[11pt]{amsart}
\usepackage[a4paper,margin=28mm]{geometry}
\usepackage{amsmath,amssymb,amsthm,amscd,mathtools,mathrsfs}
\usepackage[T1]{fontenc}
\usepackage{lmodern,microtype}
\usepackage{xcolor}
\usepackage{graphicx}
\usepackage{xurl}
\usepackage{enumitem}

\definecolor{paperlink}{RGB}{0,0,114}
\usepackage[hypertexnames=false,colorlinks=true,linkcolor=paperlink,citecolor=paperlink,urlcolor=paperlink]{hyperref}
\usepackage{lipsum}
\usepackage{adjustbox}
\usepackage{fancybox}
\usepackage{listings}
\usepackage{docmute}
\usepackage{calligra}
\usepackage{amsxtra}
\usepackage{latexsym}
\usepackage[all,cmtip]{xy}
\usepackage{comment}
\usepackage{mathabx,epsfig}
\usepackage{stmaryrd}
\usepackage{aliascnt}

\usepackage{tikz}
\usepackage{tikz-cd}
\usetikzlibrary{cd}
\usetikzlibrary{positioning, arrows.meta, matrix}
\usetikzlibrary{decorations.pathreplacing}
\tikzcdset{
  arrow style=tikz,
  arrows=-stealth
}

\makeatletter
\renewcommand\th@plain{\slshape}
\makeatother
\newtheoremstyle{plain}
 {2mm}
 {2mm}
 {\slshape}
 {}
 {\bfseries}
 {.}
 {.5em}
 {}
\theoremstyle{plain}
\newtheorem{theorem}{Theorem}[section]

\newaliascnt{corollary}{theorem}
\newtheorem{corollary}[corollary]{Corollary}
\aliascntresetthe{corollary}

\newaliascnt{lemma}{theorem}
\newtheorem{lemma}[lemma]{Lemma}
\aliascntresetthe{lemma}

\newaliascnt{proposition}{theorem}
\newtheorem{proposition}[proposition]{Proposition}
\aliascntresetthe{proposition}

\newaliascnt{claim}{theorem}
\newtheorem{claim}[claim]{Claim}
\aliascntresetthe{claim}

\newtheorem*{claim*}{Claim}

\newaliascnt{conjecture}{theorem}

\aliascntresetthe{conjecture}

\newaliascnt{question}{theorem}

\aliascntresetthe{question}

\newtheoremstyle{definition}
 {2mm}
 {2mm}
 {\normalfont}
 {}
 {\bfseries}
 {.}
 {.5em}
 {}
\theoremstyle{definition}

\newaliascnt{definition}{theorem}

\aliascntresetthe{definition}

\newaliascnt{notation}{theorem}
\newtheorem{notation}[notation]{Notation}
\aliascntresetthe{notation}

\newaliascnt{remark}{theorem}
\newtheorem{remark}[remark]{Remark}
\aliascntresetthe{remark}

\newaliascnt{example}{theorem}

\aliascntresetthe{example}

\newtheorem*{acknowledgements}{Acknowledgements}
\newtheorem*{useofai}{Use of AI tools}

\usepackage[nameinlink]{cleveref}
\usepackage{autonum}

\crefname{section}{Section}{Sections}

\crefname{theorem}{Theorem}{Theorems}
\crefname{corollary}{Corollary}{Corollaries}
\crefname{lemma}{Lemma}{Lemmas}
\crefname{proposition}{Proposition}{Propositions}
\crefname{claim}{Claim}{Claims}
\crefname{definition}{Definition}{Definitions}
\crefname{notation}{Notation}{Notations}
\crefname{problem}{Problem}{Problems}
\crefname{note}{Note}{Notes}
\crefname{remark}{Remark}{Remarks}
\crefname{example}{Example}{Examples}
\crefname{conjecture}{Conjecture}{Conjectures}
\crefname{question}{Question}{Questions}

\crefname{enumi}{}{}
\crefname{enumii}{}{}
\crefname{enumiii}{}{}
\crefformat{equation}{{\upshape(#2#1#3)}}

\numberwithin{equation}{section}

\newcommand{\PP}{\mathbb P}
\newcommand{\A}{\mathbb A}
\newcommand{\CC}{\mathbb C}
\newcommand{\OO}{\mathcal O}
\newcommand{\II}{\mathcal I}
\newcommand{\DB}{\underline{\Omega}}

\newcommand{\Supp}{\operatorname{Supp}}
\newcommand{\red}{\mathrm{red}}
\newcommand{\HH}{\mathcal H}
\newcommand{\Tor}{\operatorname{Tor}}
\newcommand{\ZZ}{\mathbb Z}
\newcommand{\Q}{\mathbb Q}

\DeclareMathOperator{\Sym}{Sym}
\DeclareMathOperator{\Frac}{Frac}
\DeclareMathOperator{\Tr}{Tr}

\DeclareMathOperator{\Ker}{Ker}
\DeclareMathOperator{\codim}{codim} 
\DeclareMathOperator{\Ima}{Im}
\DeclareMathOperator{\id}{id}
\DeclareMathOperator{\pr}{pr}
\DeclareMathOperator{\sing}{sing}
\DeclareMathOperator{\reg}{reg}
\DeclareMathOperator{\length}{length}

\newcommand{\f}{\varphi}
\renewcommand{\k}{\kappa}

\newcommand{\gp}{\mathfrak{p}}

\newcommand{\gm}{\mathfrak{m}}

\title[Totally invariant subsets of endomorphisms of projective space]{Totally invariant subsets of endomorphisms of projective space}
\author[Yohsuke Matsuzawa]{Yohsuke Matsuzawa}
\address{Department of Mathematics, Graduate School of Science, Osaka Metropolitan University, 3-3-138, Sugimoto, Sumiyoshi, Osaka, 558-8585, Japan}
\email{matsuzaway@omu.ac.jp}
\author[Shou Yoshikawa]{Shou Yoshikawa}
\address{Institute of Science Tokyo, Tokyo 152-8551, Japan}
\email{yoshikawa.s.9fe9@m.isct.ac.jp}
\date{}
\begin{document}

\begin{abstract}
We prove that every totally invariant closed subset of a non-isomorphic surjective morphism $f \colon \PP^n_{\CC} \longrightarrow \PP^n_{\CC}$ is 
a union of intersections of coordinate hyperplanes with respect to a suitable choice of homogeneous coordinates.
\end{abstract}

\maketitle

\tableofcontents

\section{Introduction}

In this paper, we prove the following theorem.

\begin{theorem}\label{thm:intro-mainthm}
    Let $f \colon \PP^n_{\CC} \longrightarrow \PP^n_{\CC}$
    be a surjective non-isomorphic morphism on complex projective space of dimension $n \geq 1$.
    Then for any Zariski closed subset $Z \subseteq \PP^n_\CC$ such that $f^{-1}(Z)=Z$ as subsets,
    there is a basis $x_0,\dots , x_n$ of $H^0(\PP^n_\CC, \OO_{\PP^n_\CC}(1))$
    such that every irreducible component $Z'$ of $Z$ is of the form
    \begin{align}
        Z'=(x_{i_1}=\cdots = x_{i_r}=0)
    \end{align}
    for some $0\leq i_1 ,\dots, i_r \leq n$.
    In other words, every irreducible component of $Z$ is an intersection of coordinate hyperplanes
    with respect to a common choice of homogeneous coordinates of $\PP^n_\CC$.
\end{theorem}

\cref{thm:intro-mainthm} immediately implies that 
there is a unique maximal proper totally invariant Zariski closed subset,
and it is a union of intersections of coordinate hyperplanes
with respect to a suitable choice of homogeneous coordinates.

We note that for a given union $Z$ of intersections of coordinate hyperplanes of $\PP^n_\CC$,
there is a non-isomorphic surjective self-morphism 
of $\PP^n$ whose maximal totally invariant closed subset
is $Z$ (\cref{sec:examples}).
Thus \cref{thm:intro-mainthm} gives a complete characterization
of the Zariski closed subsets of $\PP^n_\CC$ that are totally invariant
under some surjective non-isomorphic endomorphism.

\begin{remark}
    \cref{thm:intro-mainthm} holds over any algebraically closed field
    of characteristic zero.
    Indeed, following Lefschetz principle, 
    descent $f$ and all irreducible components of $Z$ to an algebraically 
    closed subfield $K$ of $\CC$.
    Then we can find such a homogeneous coordinates over $\CC$.
    By looking at degrees, all $Z'$ are linear over $K$ as well.
    By \cref{rem:coordlemma}, we can descent \cref{eq:linear-gluing}
    to $K$ and get a desired homogeneous coordinates over $K$.
\end{remark}

The existence of a unique maximal proper totally invariant Zariski
closed subset goes back to Dinh and Sibony's work in 2003
\cite[Section~3.4]{DinhSibony2003};
see \cite[Theorem~6.1]{DinhSibony2008} for an explicit statement.

The conjecture that every totally invariant irreducible Zariski
closed subset is linear had been established for plane curves
\cite{CerveauLinsNeto2000,FornaessSibony1994}
and for prime divisors in $\PP^3$
\cite[Corollary~1.2]{Hoering2017}
and $\PP^4$
\cite[Corollary~1.5]{BroustetHoering2026}.
It was recently proved in full generality by Luo and Meng
\cite{LuoMeng2026}.
We refer to the introduction of their paper for a more detailed
account of the history of the conjecture and earlier partial results.

Our approach, based on splitting for the relative Du Bois complex,
provides alternative proofs of the existence of maximal totally invariant subset
and the linearity result.

As an immediate corollary of \cref{thm:intro-mainthm}, we have the following.
\begin{corollary}
Let $f \colon \PP^n_{\CC} \longrightarrow \PP^n_{\CC}$ be a surjective non-isomorphic morphism on complex projective space of dimension $n \geq 1$.
Set
\[
\mathcal{Z}:=\{Z \subseteq \mathbb{P}^n_{\CC} \text{: proper closed subvariety} \mid f^{-a}(Z)=Z \text{ for some $a \in \ZZ_{>0}$}\}.
\]
Then $\# \mathcal{Z} \leq 2^{n+1}-2$.
\end{corollary}

Let us give some comments on the idea behind the proof
of \cref{thm:intro-mainthm}.
A natural way to get information about the geometric structure
of $Z$ is to look at the sheaves of differential forms $\Omega_Z^p$.
They give fairly enough information if $Z$ is smooth, for example,
but in our situation $Z$ could be singular and reducible.
A good alternative is the Du Bois complex $\DB^p_Z$,
an object of derived category of $\OO_Z$-modules with
coherent cohomologies, that encodes differential forms on smooth
models together with gluing data arising from singularities
and intersections of irreducible components.
This comes with natural restriction map
$\DB^p_{\PP^n_\CC} \longrightarrow i_* \DB^p_Z$ similar to 
the usual differential forms,
where $i \colon Z \hookrightarrow \PP^n_\CC$ is the inclusion.
The relative Du Bois complex $\DB_{\PP^n_\CC, Z}^p$ completes this morphism to
a distinguished triangle $\DB_{\PP^n_\CC, Z}^p \to \DB^p_{\PP^n_\CC} \to i_*\DB^p_Z \xrightarrow{+1}$.
When $Z$ is totally invariant, we have a strong vanishing theorem
of the sheaf cohomologies of cohomologies of the complex $\DB_{\PP^n_\CC, Z}^p$, called relative Bott vanishing.
If $Z$ is irreducible, this gives enough information 
of $\DB^p_Z$ to conclude that $Z$ is a linear subvariety.
This allows us to assume that $Z$ is a union of 
linear subvarieties. 
In this case,
$0$-th relative Du Bois complex 
$\DB^0_{\A^{n+1}, \widetilde{Z}}$,
where $\widetilde{Z}$ is the affine cone of $Z$,
has the information how $\widetilde{Z}$ is embedded in
$\A^{n+1}$. The relative Bott vanishing again imposes 
strong enough constraint to conclude $Z$ is of the form
as in \cref{thm:intro-mainthm}.

\begin{useofai}
    The proof of main theorem were generated by gpt-6-astra on October 4th, 2026.
    The example in \cref{sec:examples} was generated by ChatGPT (chat with pro mode) on November 2nd, 2026. 
    The original manuscript was written by ChatGPT, but the authors
    substantially rewrote it, as well as added explanations and lemmas.
    We take responsibility for the correctness of the contents of this manuscript.
\end{useofai}

\begin{acknowledgements}
    The first author is supported by JSPS KAKENHI Grant Numbers JP22K13903 and 23KK0252.
    The second author was supported by JSPS KAKENHI Grant number JP24K16889.
\end{acknowledgements}

\section{Du Bois complexes}
We recall the properties of Du Bois complexes used below.
Let $X$ be a reduced separated scheme of finite type over $\CC$, and put $d=\dim X$.
Let $\mathrm{Sch}/X$ denote the category of separated schemes of finite type over $X$.
The $h$-topology on this category is defined by covering families $\{T_i\to T\}$ for which $\coprod_iT_i\to T$ is a universal topological epimorphism.
Here a universal topological epimorphism is a morphism whose every base change is surjective and induces the quotient topology on the target; see~\cite[Definition~3.1]{HJ}.
In particular, Zariski open coverings and proper surjective morphisms are $h$-covers~\cite[Example~3.2]{HJ}.

Let $X_h:=(\mathrm{Sch}/X)_h$ denote the site whose underlying category is $\mathrm{Sch}/X$ and whose topology is the $h$-topology defined above.
Let $\epsilon_X:X_h\to X_{\mathrm{Zar}}$ be the canonical morphism of sites, where $X_{\mathrm{Zar}}$ denotes the small Zariski site of $X$.
For $p\geq0$, let $\Omega_{h/X}^p$ be the $h$-sheafification of the presheaf $(T\to X)\mapsto\Gamma(T,\Omega_{T}^p)$.
We define
\begin{equation}
 \DB_X^p:=R\epsilon_{X*}\Omega_{h/X}^p\,.
\end{equation}
By~\cite[Corollary~6.16]{HJ},
\begin{equation}
 \DB_X^p\in D_{\mathrm{coh}}^{[0,d]}(X),
 \qquad
 \DB_X^p\simeq0 \quad\text{for }p>d.
\end{equation}
If $X$ is smooth, then $\DB_X^p\simeq\Omega_X^p$, placed in degree zero, by~\cite[Proposition~6.14]{HJ}.

For a reduced closed subscheme $i:A\hookrightarrow X$, denote the induced morphism of $h$-sites by $i_h:A_h\to X_h$.
Following~\cite[Definition~7.24]{HJ}, put
\begin{equation}
 \Omega_{h/(X,A)}^p
 :=
 \Ker\bigl(\Omega_{h/X}^p\longrightarrow i_{h*}\Omega_{h/A}^p\bigr),
 \qquad
 \DB_{X,A}^p:=R\epsilon_{X*}\Omega_{h/(X,A)}^p\,.
\end{equation}
Here the restriction morphism $\Omega^p_{h/X} \to i_{h*}\Omega^p_{h/A}$ is induced by the pullback maps
\begin{equation}
 \Gamma(T,\Omega_{T}^p)
 \longrightarrow
 \Gamma(T\times_X A,\Omega_{(T\times_X A)}^p)
 \longrightarrow
 \Omega_{h/A}^p(T\times_X A)
\end{equation}
for objects $T\to X$ of $\mathrm{Sch}/X$.
By~\cite[Lemma~7.26]{HJ}, the restriction map is surjective and $Ri_{h*}\Omega_{h/A}^p\simeq i_{h*}\Omega_{h/A}^p$.
Applying $R\epsilon_{X*}$ to the resulting short exact sequence gives the distinguished triangle
\begin{equation}\label{eq:db-pair}
 \DB_{X,A}^p
 \longrightarrow
 \DB_X^p
 \longrightarrow
 i_*\DB_A^p
 \xrightarrow{+1}\,.
\end{equation}

Let $f:Y\to X$ be a morphism of reduced separated schemes of finite type over $\CC$, and let $f_h:Y_h\to X_h$ be the induced morphism of sites.
Let $B\subseteq Y$ and $A \subseteq X$ be reduced closed subschemes such that $f$ restricts to a morphism $f|_B \colon B\to A$.
Pullback of differential forms induces morphisms $\Omega_{h/X}^p\to f_{h*}\Omega_{h/Y}^p$ and $\Omega_{h/A}^p\to (f|_{B})_{h*}\Omega_{h/B}^p$.
Passing to derived direct images gives the natural morphisms
\begin{equation}\label{eq:db-pullback}
 \DB_X^p
 \longrightarrow
 R\epsilon_{X*}Rf_{h*}\Omega_{h/Y}^p
 \simeq
 Rf_*\DB_Y^p\, \qquad \DB_{X,A}^p
 \longrightarrow
 R\epsilon_{X*}Rf_{h*}\Omega_{h/(Y,B)}^p
 \simeq
 Rf_*\DB_{Y,B}^p\,.
\end{equation}
These constructions are compatible with identities, composition, and the distinguished triangles~\eqref{eq:db-pair}.
When $X$ and $Y$ are smooth, the first morphism in \eqref{eq:db-pullback} agrees with the morphism induced by the usual pullback of differential forms.

Let $\pi:Y\to X$ be a proper morphism with $Y$ reduced, and let $i:Z\hookrightarrow X$ be a reduced closed subscheme.
Set
\begin{equation}
 E:=(Y\times_X Z)_{\mathrm{red}},
 \qquad
 g:=\pi|_E:E\to Z.
\end{equation}
Suppose that $\pi$ induces an isomorphism $Y\setminus E\xrightarrow{\sim}X\setminus Z$.
By~\cite[Proposition~6.15]{HJ}, using invariance under reduction, there is a natural distinguished triangle
\begin{equation}\label{eq:db-blowup}
 \DB_X^p
 \longrightarrow
 R\pi_*\DB_Y^p\oplus i_*\DB_Z^p
 \longrightarrow
 i_*Rg_*\DB_E^p
 \xrightarrow{+1},
\end{equation}
whose second morphism is the difference of the natural maps.

If $X$ is integral and $\pi:Y\to X$ is a resolution of singularities, then
\begin{equation}\label{DB-canonical}
 \DB_X^{d}\xrightarrow{\simeq}R\pi_*\omega_Y.
\end{equation}
Indeed, we may choose $Z\subsetneq X$ in~\eqref{eq:db-blowup}, and the assertion follows since $\dim Z,\dim E<d$; see also~\cite[Proposition~6.19]{HJ}.

Let $X=\bigcup_{i=1}^tX_i$ be a finite covering by reduced closed subschemes.
For $\varnothing\neq I\subseteq\{1,\ldots,t\}$, set
\begin{equation}
 X_I:=\left(\bigcap_{i\in I}X_i\right)_{\mathrm{red}},
 \qquad
 i_I:X_I\hookrightarrow X,
\end{equation}
and denote the induced morphism of $h$-sites by $i_{I,h}$.
Since $\{X_i\to X\}_{i=1}^t$ is an $h$-cover, the augmented alternating \v{C}ech complex
\begin{equation}\label{eq:Cech}
 0\longrightarrow\Omega_{h/X}^p
 \longrightarrow
 \bigoplus_{|I|=1}i_{I,h*}\Omega_{h/X_I}^p
 \longrightarrow\cdots\longrightarrow
 \bigoplus_{|I|=t}i_{I,h*}\Omega_{h/X_I}^p
 \longrightarrow0
\end{equation}
is exact.
Here we use invariance of $h$-sheaves under reduction.
If every nonempty $X_I$ is smooth, applying $R\epsilon_{X*}$ gives
\begin{equation}\label{eq:db-cech}
 \DB_X^p\simeq
 \left[
 \bigoplus_{|I|=1}i_{I*}\Omega_{X_I}^p
 \longrightarrow
 \bigoplus_{|I|=2}i_{I*}\Omega_{X_I}^p
 \longrightarrow\cdots\longrightarrow
 \bigoplus_{|I|=t}i_{I*}\Omega_{X_I}^p
 \right].
\end{equation}
The summand indexed by $I$ lies in cohomological degree $|I|-1$, and the differentials are alternating restriction maps.

\begin{remark}\label{rmk:endo-union}
Let $X=\bigcup_{i=1}^tX_i$ be a finite covering by reduced closed subvarieties over $\CC$.
Assume that every nonempty $X_I$ is smooth, and let $f:X\to X$ be an endomorphism restricting to an endomorphism of each $X_i$.
Write $f_I:X_I\to X_I$ for the induced endomorphisms.
Pullback along $f$ and the $f_I$ give a commutative diagram
\begin{equation}
 \begin{tikzcd}[column sep=small]
  0 \arrow[r]
  & \Omega_{h/X}^p \arrow[r] \arrow[d]
  & \displaystyle\bigoplus_{|I|=1}i_{I,h*}\Omega_{h/X_I}^p \arrow[r] \arrow[d]
  & \cdots \arrow[r]
  & \displaystyle\bigoplus_{|I|=t}i_{I,h*}\Omega_{h/X_I}^p \arrow[r] \arrow[d]
  & 0\\
  0 \arrow[r]
  & f_{h*}\Omega_{h/X}^p \arrow[r]
  & \displaystyle\bigoplus_{|I|=1}f_{h*}i_{I,h*}\Omega_{h/X_I}^p \arrow[r]
  & \cdots \arrow[r]
  & \displaystyle\bigoplus_{|I|=t}f_{h*}i_{I,h*}\Omega_{h/X_I}^p \arrow[r]
  & 0
 \end{tikzcd}
 \,
\end{equation}
whose upper row is  exact by \eqref{eq:Cech}.
The squares commute because pullback of differential forms commutes with restriction.

Let $C^{p,\bullet}$ denote the complex in~\eqref{eq:db-cech}.
Applying $R\epsilon_{X*}$ to this diagram gives
\begin{equation}
 \begin{tikzcd}
  \DB_X^p \arrow[r, "\sim"] \arrow[d]
  & C^{p,\bullet} \arrow[d]\\
  Rf_*\DB_X^p \arrow[r, "\sim"]
  & Rf_*C^{p,\bullet}
 \end{tikzcd}
 \,.
\end{equation}
The left vertical arrow is the natural Du Bois morphism, and the right vertical arrow is induced by pullback on each $\Omega_{X_I}^p$.
If $f$ is finite, the latter is represented by the morphism of complexes $C^{p,\bullet}\to f_*C^{p,\bullet}$.
\end{remark}

\section{Linearity of totally invariant varieties}\label{sec:splitting}

In this section, we prove that every totally invariant subvariety
under a non-isomorphic surjective self-morphism of a projective space is linear.

\begin{notation}\label{not:projective}
Let $n\geq0$, and let $f:\PP^n_\CC \longrightarrow \PP^n_\CC$ be a
surjective self-morphism such that $f^*\OO_{\PP^n}(1)\simeq\OO_{\PP^n}(q)$ for an integer $q>1$.
Whenever we replace $f$ by an iterate $f^e$, we also replace $q$ by $q^e$ and retain the same notation.
\end{notation}

The following vanishing theorem is the key ingredient
of the entire proof.

\begin{proposition}[Relative Bott vanishing]\label{prop:bott}
In Notation~\ref{not:projective}, let $A\subsetneq\PP^n$ be a reduced closed subscheme with $f^{-1}(A)=A$.
For $p,j\geq0$, put $\mathcal F^{p,j}=\HH^j(\DB^p_{\PP^n,A})$.
Then for every $e\geq1$ the natural morphism
\begin{equation}
 \mathcal F^{p,j}\longrightarrow(f^e)_*\mathcal F^{p,j}
\end{equation}
is split injective, and
\begin{equation}\label{eq:bott}
 H^b(\PP^n,\mathcal F^{p,j}(t))=0\qquad(b>0,\ t>0).
\end{equation}
\end{proposition}

\begin{proof}
By~\cite[Theorem~6.9 and the proof of Corollary~6.10]{NN}, the natural morphism
\begin{equation}
 \DB_{\PP^n,A}^p\longrightarrow(f^e)_*\DB_{\PP^n,A}^p
\end{equation}
admits a splitting.
Taking cohomology, we get split injection
\begin{align}
    \mathcal{F}^{p,j} = \HH^j(\DB^p_{\PP^n,A}) \longrightarrow
    \HH^j((f^e)_* \DB^p_{\PP^n,A}) \simeq (f^e)_*\HH^j(\DB^p_{\PP^n,A}) = (f^e)_* \mathcal{F}^{p,j}.
\end{align}
Here we use $(f^e)_*$ is exact.
This proves the first statement (this morphism we just get is ``the 
natural morphism" in the statement).
Twisting by $t$ and using projection formula, 
we get split injection
\begin{align}
    \mathcal{F}^{p,j}(t) \longrightarrow
    \big((f^e)_* \mathcal{F}^{p,j}\big)(t) 
    \simeq (f^e)_* \big(\mathcal{F}^{p,j} \otimes (f^e)^*\OO_{\PP^n}(t)\big) \simeq (f^e)_* \big( \mathcal{F}^{p,j}(q^et)\big).
\end{align}
Taking sheaf cohomology, we get injection
\begin{equation}
 H^b(\PP^n,\mathcal F^{p,j}(t)) \hookrightarrow H^b(\PP^n,\mathcal F^{p,j}(q^e t)).
\end{equation}
The target vanishes for large $e$ by Serre vanishing, which finishes
the proof of the second statement.

\end{proof}

\begin{notation}\label{not:normalization}
In Notation~\ref{not:projective}, let $Z\subsetneq\PP^n$ be an integral totally invariant subvariety of dimension $d>0$, and let $i:Z\hookrightarrow\PP^n$ be the inclusion.
Let $\nu:W\to Z$ be the normalization and $\rho:\widetilde W\to W$ a resolution that is an isomorphism over the regular locus $W_{\mathrm{reg}}$ of $W$.
Write $g=f|_Z$, let $h:W\to W$ be the lift characterized by $\nu\circ h=g\circ\nu$, and put $L=\nu^*i^*\OO_{\PP^n}(1)$.
These morphisms fit into the following commutative diagram:
\begin{equation}
    \begin{tikzcd}
        & \widetilde{W} \arrow[d, "\rho"] \\
        W \arrow[r, "h"] \arrow[d, "\nu",swap]& W \arrow[d, "\nu"] \\
        Z \arrow[r, "g"] \arrow[d, hook, "i",swap]& Z \arrow[d, hook, "i"] \\
        \PP^n_\CC \arrow[r, "f"]& \PP^n_\CC
     \end{tikzcd}
\end{equation}

Since $\DB_{\PP^n}^d\simeq\Omega_{\PP^n}^d$, taking cohomology  of~\eqref{eq:db-pair} gives the exact sequence of coherent sheaves on $\PP^n$
\begin{equation}
 0\longrightarrow\HH^0(\DB^d_{\PP^n,Z})
 \longrightarrow\Omega_{\PP^n}^d
 \longrightarrow i_*\HH^0(\DB_Z^d)
 \longrightarrow\HH^1(\DB^d_{\PP^n,Z})
 \longrightarrow0.
\end{equation}
Set
\begin{alignat}{2}
 \mathcal G&:=i_*\HH^0(\DB_Z^d), \qquad&
 \mathcal M&:=\Ima(\Omega_{\PP^n}^d\to i_*\HH^0(\DB_Z^d)),\\
 \mathcal Q&:=\HH^1(\DB^d_{\PP^n,Z}), \qquad&
 \mathcal K&:=\HH^0(\DB^d_{\PP^n,Z}).
\end{alignat}

\end{notation}
Thus $\mathcal Q\simeq\mathcal G/\mathcal M$, and the exact sequence becomes
\begin{equation}\label{eq:top-sequence}
 0\longrightarrow\mathcal K\longrightarrow\Omega_{\PP^n}^d
 \longrightarrow\mathcal G\longrightarrow\mathcal Q\longrightarrow0.
\end{equation}
Taking degree-zero cohomology of \eqref{DB-canonical} yields
\begin{equation}\label{eq:top-identification}
 \mathcal G\simeq(i\circ\nu\circ\rho)_*\omega_{\widetilde W}.
\end{equation}
On $Z_{\mathrm{reg}}$, the map $\Omega^d_{\PP^n} \longrightarrow \mathcal{G}$
is the usual surjective restriction of top forms.
In particular, $\Supp\mathcal Q\subseteq Z_{\mathrm{sing}}$.

\medskip

{\bf Idea of the proof of linearity.}
We will show that $\deg Z = 1$.
To this end, consider general linear subvariety
$\Lambda \subseteq \PP^n_\CC$ of dimension $d-1$.
If we write $\overline{C} = (i \circ \nu)^{-1}(\Lambda)$,
then $\deg Z = (L \cdot \overline{C})$.
Note that if $Z$ is linear, we should be able to prove
that $\overline{C} \simeq \PP^1$ and $L|_{\overline{C}}$
is degree one line bundle.
Since $\Lambda$ is general, \cref{eq:top-sequence}
remains exact after restricted on $\Lambda$:
\begin{equation}
 0\longrightarrow\mathcal K|_\Lambda\longrightarrow\Omega_{\PP^n}^d|_\Lambda
 \longrightarrow\mathcal G|_\Lambda\longrightarrow\mathcal Q|_\Lambda\longrightarrow0.
\end{equation}
By a base change theorem and adjunction for canonical sheaves, we can show
$\mathcal G|_\Lambda \simeq 
(i \circ \nu)_* \big(\omega_{\overline{C}} \otimes L|_{\overline{C}}^{-(d-1)}\big)$.
Twisting by $d$, we have exact sequence
\begin{equation}
 0\longrightarrow\mathcal K(d)|_\Lambda\longrightarrow\Omega_{\PP^n}^d(d)|_\Lambda
 \longrightarrow
 (i \circ \nu)_* \big(\omega_{\overline{C}} \otimes L|_{\overline{C}}\big)
 \longrightarrow\mathcal Q(d)|_\Lambda\longrightarrow0.
\end{equation}
The goal is to show that 
$H^0( (i \circ \nu)_* \big(\omega_{\overline{C}} \otimes L|_{\overline{C}}\big)) 
= 
H^0(\overline{C}, \omega_{\overline{C}} \otimes L|_{\overline{C}}) = 0$, as it forces $\overline{C} \simeq \PP^1$
and $L|_{\overline{C}} \simeq \OO(1)$.
To this end, we need three inputs:
\begin{enumerate}
    \item\label{idea-Qis0} $\mathcal{Q}|_\Lambda = 0$ ;
    \item\label{idea-Pn} $H^0(\Lambda, \Omega_{\PP^n}^d(d)|_\Lambda)=0$ ;
    \item\label{idea-H1Kis0} $H^1(\Lambda, \mathcal{K}(d)|_\Lambda)=0$.
\end{enumerate}
(\cref{idea-Pn}) follows from elementary calculation of cohomology on projective space.
(\cref{idea-H1Kis0}) is a consequence of \cref{prop:bott}.
(\cref{idea-Qis0}) is also an application of \cref{prop:bott},
but it is more technical.
Note that since $\Supp \mathcal{Q} \subseteq Z_{\sing}$,
this is trivial if $Z$ is normal.
We will show that $\dim \Supp \mathcal{Q} \leq d-2$
without assuming $Z$ is normal.

\medskip

To prove $\dim \Supp \mathcal{Q} \leq d-2$,
we prepare some lemmas.

\begin{lemma}\label{lem:normal-image-is-normal}
    Consider the following diagram
    \begin{equation}
        \begin{tikzcd}
            R \arrow[r, hook, "\varphi"] \arrow[d] & S \arrow[d] \\
            R/I \arrow[r] & S/J
        \end{tikzcd}
    \end{equation}
    where $R$ is a Noetherian normal integral domain containing $\Q$,
    $S$ is an integral domain, $\varphi$ is a finite flat injective ring homomorphism,
    $I$ is an radical ideal, and $J = \sqrt{IS}$.
    If $S/J$ is normal (i.e.\ localization at every prime ideal is integrally closed integral domai), then so is $R/I$.
\end{lemma}
\begin{proof}
    Let $K = \Frac R, L = \Frac S$.
    We set $T = [L : K]^{-1}\Tr_{L/K}$, where $\Tr_{L/K}$ is the trace
    of the field extension $L/K$.
    Since $R$ is normal $\Q$-algebra and $\varphi$ is finite,
    we have $T(S) \subseteq R$.
    Note that $T \colon S \longrightarrow R$ is an $R$-module homomorphism
    with $T(1)=1$.

    We claim that $T(J) \subseteq I$.
    Indeed, let $b \in J$. Then $b^n \in IS$ for some $n \geq 1$.
    For any prime ideal $I \subseteq \mathfrak{p} \subseteq R$, 
    we have
    \begin{align}
        \Tr_{L/K}(b) = \Tr_{S_{\gp} / R_{\gp}}(b)
        \equiv \Tr_{S \otimes_R \kappa(\gp) / \k(\gp)} (b) = 0
        \mod \gp R_{\gp}
    \end{align}
    where $\k(\gp)$ is the residue field at $\gp$.
    Here for the first equality, we use $\f$ is flat so that $S_\gp$ is
    free over $R_\gp$. The last equality follows from the fact that $b$
    is nilpotent in $S \otimes_R \kappa(\gp)$.
    Therefore we get $T(b) \in \gp$, and hence $T(b) \in \sqrt{I}=I$.

    Now we get the following diagram:
    \begin{equation}
        \begin{tikzcd}
            R \arrow[r, hook, "\f"] \arrow[d] \arrow[rr, bend left = 30pt, "\id"]& S \arrow[r, "T"] \arrow[d] & R \arrow[d]\\
            R/I \arrow[r, "\overline{\f}"] \arrow[rr, bend right = 30pt, "\id",swap]& S/J \arrow[r, "\overline{T}"] & R/I
        \end{tikzcd}
    \end{equation}
    where $\overline{\f}, \overline{T}$ are induced morphisms.
    Let $\gp \subseteq R/I$ be any prime ideal.
    To end the proof, we need to show that $(R/I)_\gp$ is a normal domain.
    First note that $(S/J)_\gp$ is the product of finitely many normal local domains,
    say $(S/J)_\gp \simeq A_1 \times \cdots \times A_r$.
    Let $e_1, \dots, e_r$ be the corresponding idempotent elements of $S/J$.
    Then we have
    $\sum_{1 \leq i \leq r} \overline{T}_\gp (e_i) = \overline{T}_\gp(1) = 1$.
    Thus there is some $i$ such that $\overline{T}_\gp(e_i) \in (R/I)_\gp^\times$.
    We set
    \begin{align}
        &\psi \colon (R/I)_\gp \xlongrightarrow{\overline{\f}_\gp} (S/J)_\gp \xlongrightarrow{\pr_i} A_i\\
        &\tau \colon A_i \xlongrightarrow{\iota} (S/J)_\gp \xlongrightarrow{\overline{T}_\gp(e_i)^{-1}\overline{T}_\gp} (R/I)_\gp
    \end{align}
    where $\pr_i$ is the $i$-th projection
    and $\iota$ is the embedding to $0 \times \cdots \times A_i \times \cdots 0$.
    Then $\psi$ is a ring homomorphism, $\tau$ is an $(R/I)_\gp$-homomorphism,
    $\tau(1) = 1$. In particular, $ \tau \circ \psi = \id$.
    Since $A_i$ is integral domain, so is $(R/I)_\gp$.
    
    Let us identify $(R/I)_\gp$ as a subring of $A_i$ via $\psi$.
    Let $x = b/a \in \Frac (R/I)_\gp$ be integral over $(R/I)_\gp$,
    where $a, b \in (R/I)_\gp$ and $a \neq 0$.
    Then $x \in A_i$ since $A_i$ is normal.
    Thus $a \tau(x) = \tau(ax) = \tau(b) = b$.
    Since $\tau(x) \in (R/I)_\gp$, we get $x = b/a = \tau(x) \in (R/I)_\gp$.
\end{proof}

\begin{lemma}\label{lem:normal-image-is-normal-variety}
    Let $k$ be a field of characteristic zero.
    Consider the following diagram 
    \begin{equation}
        \begin{tikzcd}
            X \arrow[r, "f"] & Y\\
            W \arrow[u, phantom, "\subseteq", sloped] \arrow[r, "g"] & Z \arrow[u, phantom, "\subseteq", sloped]
        \end{tikzcd}
    \end{equation}
    where $X, Y$ are normal schemes of finite type over $k$,
    $f$ is a finite flat morphism, $Z \subseteq Y$ is a reduced closed subscheme,
    and $W = f^{-1}(Z)_{\red}$.
    Let $x \in W$ be a scheme point and $y= g(x)$. 
    If $\OO_{W, x}$ is normal, then so is $\OO_{Z, y}$.
\end{lemma}
\begin{proof}
    We have the following diagram of completed local rings
    \begin{equation}
        \begin{tikzcd}
            \widehat{\OO}_{Y, y} \arrow[r, "f^*"] \arrow[d] & \widehat{\OO}_{X, x} \arrow[d] \\
            \widehat{\OO}_{Z, y} \arrow[r,"g^*"] & \widehat{\OO}_{W, x}
        \end{tikzcd}
    \end{equation}
    where the vertical arrows are surjective.
    Since the local rings of schemes of finite type over a field are
    excellent (and in particular Nagata), $\widehat{\OO}_{X, x}, \widehat{\OO}_{Y, y}, \widehat{\OO}_{W, x}$
    are normal \cite[\href{https://stacks.math.columbia.edu/tag/0C23}{Tag 0C23}]{stacks-project},
    and $\widehat{\OO}_{Z, y}$ is reduced \cite[\href{https://stacks.math.columbia.edu/tag/07NZ}{Tag 07NZ}]{stacks-project}. 
    Moreover, since $f$ is flat finite, $f^*$ is a flat finite ring homomorphism
    (and in particular injective).
    By \cref{lem:normal-image-is-normal}, $\widehat{\OO}_{Z, y}$ is normal.
    By the faithfully flat descent of normality, $\OO_{Z, y}$ is also normal \cite[\href{https://stacks.math.columbia.edu/tag/0FIZ}{Tag 0FIZ}]{stacks-project}.
\end{proof}

\begin{lemma}\label{lem:branches}
In Notation~\ref{not:normalization}, let $N\subseteq Z$ be the nonnormal locus, and let $B\subseteq W$ be the union of the codimension-one components of $\nu^{-1}N$.
Then $g^{-1}N = N$ and $h^{-1}(B)=B$.
After replacing $f$ with some iterate, 
each irreducible component $E$ of $B$ satisfies $h^*E=qE$.
\end{lemma}

\begin{proof}
By \cref{lem:normal-image-is-normal-variety} or \cite[Theorem 1.2]{Lebl}, we have $g^{-1}N\subseteq N$.
Since $g$ is finite, we have $g^{-1}(N)=N$.
Consequently, $h^{-1}(B)= B$.
Thus $h^{-1}$ permutes the components of $B$.
Replacing $f$ with some iterate, we may assume 
each component of $B$ is totally invariant under $h$.

Let $E \subseteq B$ be an irreducible component of $B$.
Then we can write $h^*E=aE$ for some $a \geq 1$.
We have
\begin{equation}
 \deg(h)(L^{d-1}\cdot E) =((h^*L)^{d-1}\cdot (h^*E))=q^{d-1}(L^{d-1}\cdot h^*E) = q^{d-1}a(L^{d-1}\cdot E),
\end{equation}
where the first equality holds since $L$ is base point free.
As $(L^{d-1} \cdot E) > 0$ and 
$\deg(h)=q^d$, we have $h^*E=qE$, as desired.
\end{proof}

\begin{proposition}\label{prop:codim-two}
In Notation~\ref{not:normalization}, one has $\dim\Supp\mathcal Q\leq d-2$.
\end{proposition}

\begin{proof}
Let $N$ and $B$ be as in Lemma~\ref{lem:branches}, and replace $f$ by the iterate provided there.
Let $U$ be the complement of the union of $\nu(W_{\sing})$ and irreducible components of $N$ of codimension at least $2$ in $Z$.
Then $U \cap N$ is pure codimension one in $U$, the scheme $\nu^{-1}(U)$ is regular, and $U$ contains all codimension one points of $Z$.
Since $\nu^{-1}(U) \to U$ is a resolution and the normalization, we have
\begin{equation}
N \cap U=Z_{\sing} \cap U.
\end{equation}
Thus we have
\[
\Supp(\mathcal Q|_{U}) \subseteq Z_{\sing} \cap U=N \cap U.
\]
Then there exists a positive integer $c$ such that $\mathcal I_N^c \mathcal Q|_U =0$, where $\mathcal I_N$ is the defining ideal of $N$.
By the exact sequence \eqref{eq:top-sequence}, we have
\begin{equation}\label{eq:incl1}
\mathcal I_N^c \mathcal H^0(\DB^d_Z)|_U \subseteq \mathcal M|_U.
\end{equation}
Let $N' \subseteq N$ be the union of all codimension one irreducible component of $N$.
Set $B:=\nu^{-1}(N')_{\red}$ as the reduced divisor on $W$.
Then $B|_{\nu^{-1}(U)}=\nu^{-1}(U \cap N)$ by the choice of $U$.
Then there exists a positive integer $m$ such that
\begin{equation}\label{eq:incl2}
\OO_{\nu^{-1}(U)}(-(q^m-1)B) \subseteq \II^c_N \OO_{\nu^{-1}(U)}.
\end{equation}
Therefore, we obtain
\begin{equation}\label{eq:incl3}
\nu_*\omega_{\nu^{-1}(U)}(-(q^m-1)B) \subseteq \nu_*(\II^c_N\omega_{\nu^{-1}(U)}) =\II_N^c \mathcal H^0(\DB^d_Z)|_U \subseteq \mathcal M|U,
\end{equation}
where the first inclusion follows from \eqref{eq:incl2}, the equation follows from \eqref{eq:top-identification} and the smoothness of $\nu^{-1}(U)$, and the second inclusion follows from \eqref{eq:incl1}.

It follows from \Cref{lem:branches} that $h^*B=qB$.
In particular, we have $R_{h^m} \geq (q^m-1)B$, where $R_{h^m}$ is the ramification divisor of $h^m$.
Therefore, the natural map induces
\[
\omega_W \to (h^m)_*(\omega_W(-R_{h^m})) \hookrightarrow (h^m)_*(\omega_W(-(q^m-1)B)).
\]
Set 
$U' := U \cap (Z \setminus g^m(Z \setminus U))$, then
\[
g^{-m}(U') \subseteq U
\]
and $U'$ contains all codimension one points.
Then the natural map induces
\begin{equation}\label{eq:incl4}
\nu_*\omega_{\nu^{-1}(U')} \to \nu_*h^m_*\omega_W(-(q^m-1)B)|_{U'}=g^m_*\bigl(\nu_*(\omega_W(-(q^m-1)B)|_{g^{-m}(U')})\bigr) \subseteq g^m_*(\mathcal M|_{g^{-m}(U')}),
\end{equation}
where the inclusion follows from \eqref{eq:incl3}.
Since $\nu^{-1}(U')$ is regular, we have $\mathcal H^0(\DB^d_Z)|_{U'}=\nu_*\omega_{\nu^{-1}(U')}$.
Considering the commutative diagram in which each horizontal sequence is exact;
\begin{equation}
\begin{tikzcd}
0 \arrow[r] & \mathcal{M}|_{U'} \arrow[r] \arrow[d] & \mathcal{H}^0(\DB^d_Z)|_{U'} \arrow[r] \arrow[d] & \mathcal{Q}|_{U'} \arrow[r] \arrow[d] & 0 \\
0 \arrow[r] & (g^m_*\mathcal{M})|_{U'} \arrow[r] &  (g^m_*\mathcal{H}^0(\DB^d_Z))|_{U'} \arrow[r]  & (g^m_*\mathcal{Q})|_{U'} \arrow[r]  & 0.
\end{tikzcd}
\end{equation}
By \eqref{eq:incl4} and $\mathcal H^0(\DB^d_Z)|_{U'}=\nu_*\omega_{\nu^{-1}(U')}$, the middle vertical map factors through $(g^m_*\mathcal{M})|_{U'}$.
In particular, the right vertical map is zero.
Applying \cref{prop:bott} to $\mathcal{Q}=\mathcal{H}^1(\DB^d_{\PP^n,Z})$, we obtain $\mathcal{Q}|_{U'}=0$.
Since $U'$ contains all codimension one points of $Z$, we obtain $\dim \Supp \mathcal{Q} \leq d-2$. 
\end{proof}

\begin{theorem}\label{thm:linear}
In Notation~\ref{not:normalization}, $Z$ has degree one and is a projective linear subspace.
\end{theorem}

\begin{proof}
Set $\pi:=i\circ\nu:W\to\PP^n$.
Choose a general linear subspace $\Lambda\simeq\PP^{n-d+1}$ of codimension $r:=d-1$.
By Proposition~\ref{prop:codim-two} and the normality of $W$, both $\Supp\mathcal Q$ and $\pi(W_{\mathrm{sing}})$ have dimension at most $d-2$.
Thus $\Lambda$ avoids both closed subsets.

Let $\overline C=\pi^{-1}(\Lambda)$, the scheme-theoretic inverse image.
For $d>1$, this is the complete intersection of $r$ general members of the basepoint-free linear system on $W$ induced by $|\OO_{\PP^n}(1)|$.
Bertini's theorem shows that $\overline C$ is a smooth integral curve contained in $W_{\mathrm{reg}}$.
For $d=1$, we have $\Lambda=\PP^n$ and $\overline C=W$, which is already a smooth integral curve.
Let $\alpha:\overline C\to\Lambda$ be the finite morphism induced by $\pi$.
Since $\overline C$ is cut out by $d-1$ divisors associated with $L$ and $\nu$ is birational, we have
\begin{equation}\label{eq:curve-degree}
 \deg(L|_{\overline C})=L^d=\deg Z.
\end{equation}

Since $\Lambda$ is general, the Kleiman--Bertini theorem for sheaf tensor products~\cite[Theorem]{MillerSpeyer} gives
\begin{equation}\label{eq:tor-van}
 \Tor_j^{\OO_{\PP^n}}(\mathcal F,\OO_\Lambda)=0
 \qquad(j>0,\ \mathcal F\in\{\mathcal G,\mathcal K\})\,.
\end{equation}
 Take an open subset $U$ of $\PP^n$ containing $\Lambda$ such that $\overline{C} \subseteq \pi^{-1}(U)\subseteq W_{\reg}$.
Applying \cite[\href{https://stacks.math.columbia.edu/tag/02KG}{Tag 02KG}]{stacks-project} to the pullback diagram
\begin{equation}
\begin{tikzcd}
\overline{C} \arrow[r,"\alpha"] \arrow[d] & \Lambda \arrow[d] \\
\pi^{-1}(U) \arrow[r,"\pi"] & U,
\end{tikzcd}
\end{equation}
we obtain
\begin{equation}
 \mathcal G|_\Lambda \simeq (\pi_*\omega_{\pi^{-1}(U)})|_{\Lambda}
 \simeq\alpha_*(\omega_W|_{\overline C})
 \simeq\alpha_*(\omega_{\overline C}\otimes L^{-r}|_{\overline C}),
 \label{eq:restGtoLambda}
\end{equation}
where the last equivalence follows from the adjunction on $W_{\mathrm{reg}}$.

We restrict \cref{eq:top-sequence} to $\Lambda$.
By $\Lambda  \cap \Supp\mathcal Q = \emptyset$, \cref{eq:tor-van}, and
\cref{eq:restGtoLambda},
after twisting by $d$ and using projection formula, we get
\begin{equation}\label{eq:curve-sequence}
 0\longrightarrow(\mathcal K|_\Lambda)(d)
 \longrightarrow\Omega_{\PP^n}^d(d)|_\Lambda
 \longrightarrow\alpha_*(\omega_{\overline C}\otimes L|_{\overline C})
 \longrightarrow0\,.
\end{equation}

Now, we prove the following claim.
\begin{claim}\label{claim:vanishing}
We have
\[
H^1(\Lambda,\mathcal{K}(d)|_{\Lambda})=H^0(\Omega^d_{\PP^n}(d)|_{\Lambda})=0.
\]
\end{claim}
\begin{proof}[Proof of \Cref{claim:vanishing}]
Consider a coherent sheaf $\mathcal{F}$ such that $\Tor_j^{\OO_{\PP^n}}(\mathcal{F},\OO_{\Lambda})=0$ for $j>0$.
We consider the Koszul resolution 
\[
\OO_{\PP^n}(-(d-1))\to \cdots \to \mathcal{O}_{\PP^n}(-1)^{\oplus d-1} \to \OO_{\PP^n}
\]
of $\OO_{\Lambda}$.
Tensoring the complex with $\mathcal{F}$ over $\OO_{\PP^n}$, we obtain the Koszul resolution of $\mathcal{F}|_{\Lambda}$ by the Tor vanishing.
Thus, we obtain implications
\begin{align}
H^{a+1}(\mathcal{F}(-a))=0 \quad (d-1 \geq a \geq 0) &\Rightarrow H^1(\mathcal{F}|_{\Lambda})=0 \\
H^{a}(\mathcal{F}(-a))=0 \quad (d-1 \geq a \geq 0) &\Rightarrow H^0(\mathcal{F}|_{\Lambda})=0.
\end{align}
Thus, the vanishing of $H^1(\Lambda,\mathcal{K}(d)|_{\Lambda})$ follows from \eqref{eq:tor-van} and \cref{prop:bott}.
The vanishing of $H^0(\Omega^d_{\PP^n}(d)|_{\Lambda})$ follows from the flatness of $\Omega^d_{\PP^n}(d)$ and \Cref{appendix-prop:bottvanforprojsp}.
\end{proof}

The cohomology sequence of~\eqref{eq:curve-sequence} consequently gives $H^0(\overline C,\omega_{\overline C}\otimes L|_{\overline C})=0$.
Since $\deg(L|_{\overline C})>0$, 
$\overline C$ is a projective line.
Thus, $\deg Z=1$.
Hence $Z\subseteq\PP^n$ is a linear subspace by~\cite[Chapter~I, Exercise~7.6]{Hartshorne} or \cref{lem:degoneislin}.
\end{proof}

\section{Totally invariant sets are a union of intersections of coordinate hyperplanes}\label{sec:graded}

In this section, we prove that if a finite union of linear subvarieties
is totally invariant under a non-isomorphic surjective self-morphism of a projective 
space, these linear subvarieties are intersection of coordinate 
hyperplanes with respect to a common homogeneous coordinates.
Combined with \cref{thm:linear}, 
this complete the proof of \cref{thm:intro-mainthm}.

\begin{notation}\label{not:affine}
Let $V$ be a finite-dimensional $\CC$-vector space, and set 
\begin{align}
    S=\operatorname{Sym}_{\CC}V^* = \bigoplus_{m \geq 0} \Sym_\CC^m V^*
\end{align}
with the natural graded ring structure.
Let $\mathfrak m=S_{>0}$ be the irrelevant ideal.
For a vector subspace $U\subseteq V$, write $\mathbb A(U)=\operatorname{Spec}(\operatorname{Sym}_{\CC}U^*)$ and $\PP(U)=\operatorname{Proj}(\operatorname{Sym}_{\CC}U^*)$.
Then the quotient $S\twoheadrightarrow\operatorname{Sym}_{\CC}U^*$ defines the closed affine subspace $\mathbb A(U)\subseteq\mathbb A(V)$ and the projective linear subspace $\PP(U)\subseteq\PP(V)$.
For a graded $S$-module $M$, the degree $l$ part will be 
denoted by $M_l$.
\end{notation}

We first prove a purely linear-algebraic lemma that provides a sufficient condition for a family of linear subspaces to be simultaneously spanned by subsets of a common basis.

\begin{lemma}\label{lem:coordinates}
Let $V$ be a finite-dimensional $\CC$-vector space, and let $\mathscr U$ be a finite collection of vector subspaces of $V$, closed under intersections and containing $0$.
Suppose for every nonempty subfamily $\mathscr V$,
\begin{equation}\label{eq:linear-gluing}
 V^* \longrightarrow
 \bigoplus_{U \in \mathscr V}U^*\longrightarrow
 \bigoplus_{U, U' \in \mathscr V}(U\cap U')^*
\end{equation}
is exact, where the first map is defined by restriction,
and the second is differences of restrictions.
Then there is a basis of $V$ in which every member of $\mathscr U$ is spanned by a subset of that basis.
\end{lemma}

\begin{proof}
Dualizing~\eqref{eq:linear-gluing} gives a exact sequence
\begin{equation}
 \bigoplus_{U,U' \in \mathscr V}(U\cap U')
 \xrightarrow{\partial}\bigoplus_{U \in \mathscr V}U
 \xrightarrow{\sum}V\,,
\end{equation}
where $\Sigma((u_U)_U) = \sum_{U}u_U$,
and $\partial$ is such that for $u \in U \cap U'$,
\begin{align}
    \pr_U(\partial(u)) = -u,\quad \pr_{U'}(\partial(u))=-u,\quad
    \pr_{U''}(\partial(u))=0 \ \text{for $U'' \neq U, U'$}.
\end{align}
Applying this to $\mathscr{V}=\{U,U_1,\ldots,U_t\}\subseteq\mathscr U$, we get
\begin{equation}\label{eq:distribute}
 U\cap\sum_iU_i=\pr_U (\Ker \Sigma) = \pr_U(\Ima \partial)  =\sum_i(U\cap U_i)\,.
\end{equation}

For each $U\in\mathscr U$, choose $W_U$ such that
\begin{equation}
 U=W_U\oplus
 \sum_{\substack{U'\in\mathscr U\\U'\subsetneq U}}U'\,.
\end{equation}
Induction on inclusion gives
\begin{equation}
 U=\sum_{\substack{U'\in\mathscr U\\U'\subseteq U}}W_{U'}\,.
\end{equation}
We show that the nonzero subspaces among $W_U$ are independent.
Indeed, suppose that
\begin{equation}
 w_1+\cdots+w_s=0,
 \qquad
 0\ne w_i\in W_{U_i},
\end{equation}
where the $U_i$ are distinct, and choose $U_1$ maximal among them under inclusion.
By~\eqref{eq:distribute},
\begin{equation}
 w_1\in U_1\cap\sum_{i=2}^sU_i
 =\sum_{i=2}^s(U_1\cap U_i)\,.
\end{equation}
Since $\mathscr U$ is closed under intersections and $U_1$ is maximal, each $U_1\cap U_i$ belongs to $\mathscr U$ and is strictly contained in $U_1$.
Hence
\begin{equation}
 w_1\in
 W_{U_1}\cap
 \sum_{\substack{U'\in\mathscr U\\U'\subsetneq U_1}}U'
 =0,
\end{equation}
a contradiction.
Thus
\begin{equation}
 U=\bigoplus_{\substack{U'\in\mathscr U\\U'\subseteq U}}W_{U'}
 \qquad(U\in\mathscr U)\,.
\end{equation}
Choose bases of all $W_U$ and extend their union to a basis of $V$.
Every member of $\mathscr U$ is spanned by a subset of this basis.
\end{proof}

\begin{remark}\label{rem:coordlemma}
    Converse of \cref{lem:coordinates} also holds.
    Also, it is obviously true over any fields.
\end{remark}

We will apply \cref{lem:coordinates} to a family of
affine cone of totally invariant linear subvarieties.

Let $f \colon \PP(V) \longrightarrow \PP(V)$ be a surjective morphism 
with $f^*\OO(1) \simeq \OO(q)$ for some integer $q>1$.
Let $F \colon \A(V) \to \A(V)$ be its homogeneous lift, so that
$F$ is a finite surjective morphism such that
$F^{\#}(S_1) \subseteq S_q$.
Let $U_1,\ldots,U_t\subseteq V$ be $\CC$-vector subspaces such that 
\begin{align}
    F^{-1}(\A(U_i)) \subseteq \A(U_i).
\end{align} 
Put $A=\bigcup_i\mathbb A(U_i)$, endowed with its reduced closed subscheme structure.

A key observation, which we explain below, is that the sequence \cref{eq:linear-gluing}
with $\mathscr{V}=\{U_1, \dots, U_t\}$
shows up in a degree one part of relative Du Bois complex 
$\DB^0_{\A(V),A}$.
Therefore, our task is to show the vanishing of cohomology of 
the relative Du Bois complex.

For a vector subspace $U\subseteq V$, we have
\begin{equation}
 \Omega^p_{\mathbb A(U)}
 \simeq
 \operatorname{Sym}_{\CC}U^*\otimes_{\CC}\bigwedge^pU^*.
\end{equation}
This has natural graded $\Sym_\CC U^*$-module structure with
$\deg(du)=1$ for $u\in U^*$.
In particular, we have $\Omega^p_{\A(U)} \simeq (\Sym_\CC U^*)(-p)^{\oplus \binom{\dim U}{p}}$ as graded modules.

By \eqref{eq:db-cech}, we have
\begin{equation}\label{eq:cech}
 \DB^p_{\mathbb A(V),A}
 \simeq
 \left[
 \Omega^p_{\mathbb A(V)}
 \longrightarrow
 \bigoplus_i\Omega^p_{\mathbb A(U_i)}
 \longrightarrow
 \bigoplus_{i<j}\Omega^p_{\mathbb A(U_i\cap U_j)}
 \longrightarrow\cdots
 \right]\,,
\end{equation}
where the differentials are alternating restriction maps.
We endow the graded structure on $\mathcal H^i(\DB^p_{\mathbb A(V),A})$ induced by the complex \eqref{eq:cech} of graded $S$-modules.

For $p=0$, the first three terms of \cref{eq:cech} are
\begin{align}
    S = \Sym_\CC V^* \longrightarrow \bigoplus_i \Sym_\CC U_i^*
    \longrightarrow \bigoplus_{i<j} \Sym_\CC (U_i \cap U_j)^*.
\end{align}
Its degree one part is exactly \cref{eq:linear-gluing},
up to the choice of indexing by strictly increasing tuples or all tuples.
We will show that the first cohomology
\begin{equation}
 \HH^1(\DB^0_{\mathbb A(V),A})
 =
 \frac{
 \left\{(g_i)\in\bigoplus_i\operatorname{Sym}_{\CC}U_i^*
 \;\middle|\;
 g_i|_{U_i\cap U_j}=g_j|_{U_i\cap U_j}\right\}
 }{
 \left\{(g|_{U_i})_i\mid g\in S\right\}
 }\label{eq:H1Omega0}
\end{equation}
actually vanishes.
We note that we automatically have
\begin{align}
    \HH^1(\DB^0_{\mathbb A(V),A})_\ell=0 \quad \text{for $\ell \leq 0$.}
    \label{eq:van-H1-nonpos-deg}
\end{align}
Indeed, every $U_i$ contains the origin, so compatible constant functions extend to a constant function on $\mathbb A(V)$. This implies $ \HH^1(\DB^0_{\mathbb A(V),A})_0=0$.
Since the complex \cref{eq:cech} has no negative graded pieces, 
we obtain $ \HH^1(\DB^0_{\mathbb A(V),A})_\ell=0$ for $\ell < 0$.


Recall that the local cohomology module
$H^a_{\mathfrak m}(\mathcal H^i(\DB^j_{\A(V),A}))$
is naturally graded over $S$ (see for example \cite[p.~187]{GW} or \cite[Chapter 13]{Brodmann_Sharp_2012}).

\begin{lemma}\label{lem:positive-weights}
In the above notation, we have
\begin{equation}\label{eq:positive-local-cohomology}
 H^a_{\mathfrak m}(\mathcal H^i(\DB^j_{\A(V),A}))_\ell=0\qquad(a,i,j\geq0,\ \ell>0)\,.
\end{equation}
\end{lemma}

\begin{proof}
The affine space $\A(V)$ has the natural $\mathbb{G}_m$-action and $F$ is $\mathbb{G}_m$-equivariant by the condition $F^{\#}(S_1) \subseteq S_q$.
Thus, the fiber $F^{-1}(0)$ is finite and stable under scalar multiplication, so it consists of the origin.
Thus $\sqrt{(F^e)^\#(\mathfrak m)S}=\mathfrak m$.
By the proof of~\cite[Corollary 6.10]{NN}, the natural map 
\begin{align}
    \DB^j_{\A(V),A} \longrightarrow F^e_*\DB^j_{\A(V),A}
    \label{eq:splitinjDBFepush}
\end{align}
splits.
Here the right hand side is simply $\DB^j_{\A(V),A}$,
but each term has the $S$-module structure induced by $S \xrightarrow{(F^e)^\#} S$.
Furthermore, by \Cref{rmk:endo-union}, the above natural map commutes with the morphism of complexes induced by the natural map
\[
\Omega^j_{\A(U_{i_1}\cap \cdots \cap U_{i_l})} \longrightarrow F^e_*\Omega^j_{\A(U_{i_1}\cap \cdots \cap U_{i_l})},
\]
which sends degree $t$ elements to degree $q^et$ elements.
Since the graded structure on $\DB^j_{\A(V), A}$ is 
defined by \cref{eq:cech}, \cref{eq:splitinjDBFepush}
sends degree $t$ elements to degree $q^et$ elements as well.
Taking complex cohomology of \cref{eq:splitinjDBFepush},
we get split injection
\begin{align}
    \HH^i(\DB^j_{\A(V),A}) \longrightarrow
    F^e_*\HH^j(\DB^j_{\A(V),A})
\end{align}
sending degree $t$ elements to degree $q^et$ elements.
Since $\sqrt{(F^e)^\#(\mathfrak m)S}=\mathfrak m$,
we have
\begin{align}
    H^a_{\gm}( F^e_*\HH^j(\DB^j_{\A(V),A})) \simeq
    F^e_* H^a_{\gm}(\HH^j(\DB^j_{\A(V),A})) 
\end{align}
as graded $S$-module for all $a \geq 0$.
Thus, we obtain an injection
\begin{equation}
 H^a_{\mathfrak m}(\mathcal H^i(\DB^j_{\A(V),A}))_\ell\hookrightarrow H^a_{\mathfrak m}(\mathcal H^i(\DB^j_{\A(V),A}))_{q^e\ell}\,.
\end{equation}
For each $a$, the graded module $H^a_{\mathfrak m}(\mathcal H^i(\DB^j_{\A(V),A}))$ is bounded above in degree since it is an Artinian $S$-module.
For $\ell>0$, the target therefore vanishes for sufficiently large $e$.
Thus, we have $H^a_{\mathfrak m}(\mathcal H^i(\DB^j_{\A(V),A}))_\ell=0$, as desired.
\end{proof}

\begin{proposition}\label{prop:arrangement-db}
In the above notation, we have
\begin{equation}\label{eq:arrangement-db}
 \HH^1(\DB^0_{\mathbb A(V),A})=0\,.
\end{equation}
\end{proposition}

\begin{proof}
Put $D_A:=\HH^1(\DB^0_{\mathbb A(V),A})$.
We prove by induction on $t$.
For $t=1$, the restriction $S\to\operatorname{Sym}_{\CC}U_1^*$ is surjective, so $D_A=0$.
Suppose $t>1$ and put $C:=\bigcap_iU_i$ and $c:=\dim C$.
Take $x \in \A(V) \backslash \A(C)$, then there exists $j$ such that $x \notin \A(U_j)$ since $\A(C)=\A(U_1) \cap \cdots \cap \A(U_t)$.
By \eqref{eq:cech}, we have
\[
\HH^1(\DB^0_{\mathbb A(V),A})_x \simeq \HH^1(\DB^0_{\mathbb A(V),\bigcup_{i \neq j}\A(U_i)})_x,
\]
which vanishes by the induction hypothesis.
Thus, we have
\begin{equation}\label{eq:supp}
 \Supp D_A\subseteq\mathbb A(C)\,.
\end{equation}

Choose a $\CC$-vector subspace $V' \subseteq V$ with $V=C \oplus V'$, and put 
$S_C:=\operatorname{Sym}_{\CC}C^*$ and $S':=\operatorname{Sym}_{\CC}(V')^*$
so that $S=S_C\otimes_{\CC}S'$.
Let $U_i':=U_i\cap V'$ and $A':=\bigcup_i\mathbb A(U_i')$.
Since $U_i=C\oplus U_i'$ for every $i$, we have $\mathbb A(U_i)=\mathbb A(C)\times\mathbb A(U_i')$.
Taking the union over $i$, we obtain
\begin{equation}
 A=\mathbb A(C)\times A'.
\end{equation}
The equality holds as closed subschemes since $\mathbb A(C)\times A'$ is reduced, which is because the coordinate ring is 
a polynomial ring over a reduced ring.

Applying this product decomposition termwise to~\eqref{eq:cech} gives
\begin{equation}\label{eq:Kunneth}
 \DB^p_{\mathbb A(V),A}
 \simeq
 \bigoplus_{a+b=p}
 \Omega^a_{\mathbb A(C)}
 \otimes_{\CC}\DB^b_{\mathbb A(V'),A'}\,.
\end{equation}
Setting $p=0$ and taking $\mathcal{H}^1$, we obtain
\begin{equation}\label{eq:product}
 D_A\simeq S_C\otimes_{\CC}D_{A'},
 \qquad
 D_{A'}:=\HH^1(\DB^0_{\mathbb A(V'),A'})\,.
\end{equation}
Thus, $\Supp D_A = \mathbb{A}(C) \times \Supp D_{A'}$ and
\eqref{eq:supp} show that the finite graded $S'$-module $D_{A'}$ is supported at the origin, hence has finite length.
Moreover, \cref{eq:van-H1-nonpos-deg} applies to
$D_{A'}$ as well, and thus
$(D_{A'})_\ell=0$ for $\ell\leq0$.

By \eqref{eq:Kunneth} for $p=c$, we have the splitting injection
\[
\omega_{\mathbb A(C)}\otimes_{\CC}D_{A'} \hookrightarrow P:=\HH^1(\DB^c_{\mathbb A(V),A}).
\]
By Lemma~\ref{lem:positive-weights}, we have
\begin{equation}\label{eq:van-l}
 H^c_{\mathfrak m}
 \bigl(\omega_{\mathbb A(C)}\otimes_{\CC}D_{A'}\bigr)_\ell=0
 \qquad(\ell>0)\,.
\end{equation}

Let $z_1,\ldots,z_c$ be linear coordinates on $C$, and put $\mathfrak m_C=(z_1,\ldots,z_c)\subseteq S_C$.
Since $D_{A'}$ has finite length,
we have $H^i_{S'_{>0}}(D_{A'})=0$ for $i \neq 0$ and 
$H^0_{S'_{>0}}(D_{A'})=D_{A'}$.
Thus the K\"unneth formula for local cohomology~\cite[Theorem~(2.2.5)]{GW} gives an isomorphism of graded $S$-modules
\begin{equation}
 H^c_{\mathfrak m}
 \bigl(\omega_{\mathbb A(C)}\otimes_{\CC}D_{A'}\bigr)
 \simeq
 H^c_{\mathfrak m_C}(\omega_{\mathbb A(C)})
 \otimes_{\CC}D_{A'}.
\end{equation}
In particular, we have $H^c_{\mathfrak m_C}(\omega_{\mathbb A(C)})_0
 \otimes_{\CC}(D_{A'})_{\ell}=0$ for $\ell>0$ by \eqref{eq:van-l}.
Since the class
\begin{equation}
 \theta=
 \frac{dz_1\wedge\cdots\wedge dz_c}{z_1\cdots z_c}
 \in H^c_{\mathfrak m_C}(\omega_{\mathbb A(C)})
\end{equation}
is nonzero and has degree zero, we have $(D_{A'})_{\ell}=0$ for $\ell>0$.
Combining $(D_{A'})_{\ell}=0$ for $\ell \leq 0$, we have $D_{A'}=0$.
It follows from \eqref{eq:product} that $D_A=0$.
\end{proof}

\begin{proposition}\label{prop:linear-configuration}
Let $f \colon \PP^n_\CC \longrightarrow \PP^n_\CC$ be 
a non-isomorphic surjective self-morphism.
Let $Z_1,\ldots,Z_t \subseteq \PP^n_\CC$ be linear subvarieties 
such that $f^{-1}(Z_i)=Z_i$ for $i=1,\dots, t$.
Then there exists a coordinate $x_0,\ldots,x_n$ on $\PP^n_\CC$ such that $Z_i$ is an intersection of some members in $V_+(x_0),\ldots,V_+(x_n)$ for every $i$.
\end{proposition}

\begin{proof}
Let $V=\CC^{n+1}$, identify $\PP^n$ with $\PP(V)$.
The endomorphism $f$ induces the endomoprhism $F \colon \A(V) \to \A(V)$ such that $F^{\#}(V^*) \subseteq \Sym_\CC^q V^*$.
Let $U_i$ be a vector subspace of $V$ such that $\PP(U_i)=Z_i$ for $i=1,\ldots,t$.
Total invariance of $Z_i$ gives $F^{-1}(\mathbb A(U_i))=\mathbb A(U_i)$.

We apply \cref{lem:coordinates} to 
$V$ and $\mathscr{U} = \{\text{intersections of $U_i$}\} \cup \{0\}$. 
Note that all the members of $\mathscr{U}$ are totally invariant
under $F$.
The assumption of \cref{lem:coordinates} is 
rephrased as the vanishing of degree one part of
\cref{eq:H1Omega0} for $A = \bigcup_{U \in \mathscr{V}} \A(U)$
where $\mathscr{V}$ is any subfamily of $\mathscr{U}$.
\cref{prop:arrangement-db} guaranties this vanishing.
Thus \cref{lem:coordinates} gives a basis $v_0,\ldots,v_n$ of $V$ such that $U_1,\ldots,U_t$ are spanned by the basis. Let $x_0,\dots, x_n \in V^*$ be the dual basis. 
Then each $Z_i$ is an intersection of hyperplanes defined by members of $x_0,\ldots,x_n$.
\end{proof}

\begin{proof}[Proof of Theorem~\ref{thm:intro-mainthm}]
Replace $f$ by an iterate such that each irreducible component of $Z$ is totally invariant.
By \cref{thm:linear}, every irreducible component is linear, and \cref{prop:linear-configuration} gives desired a coordinate $x_0,\ldots,x_n$. 
\end{proof}

\section{Example}\label{sec:examples}

For a given union of intersections of coordinate hyperplanes $Z$
of $\PP^n_\CC$, 
there exist a non-isomorphic surjective self-morphism 
of $\PP^n$ whose maximal totally invariant closed subset
is $Z$.
To construct examples, let us use the following natation.
For a subset $A \subseteq \{0,\dots, n\}$, we set
\begin{align}
    L_A = \bigcap_{i \in A}V_+(x_i) \subseteq \PP^n_\CC
\end{align}
where $V_+(x_i)$ is the closed subscheme defined by $x_i$.
Subvariety of the form $L_A$ is called coordinate linear subvariety.

Let us fix a finite family $\{A_\lambda\}_{\lambda \in \Lambda
}$ of subsets of $\{0,\dots, n\}$ such that 
$A_\lambda \neq \emptyset, \{0,\dots, n\}$.
Then we set
\begin{align}
    Z := \bigcup_{\lambda \in \Lambda} L_{A_\lambda},\quad
    I = \bigcap_{\lambda \in \Lambda} (x_j \mid j \in A_\lambda)
\end{align}
where $I$ is a homogeneous ideal of $\CC[x_0,\dots, x_n]$ defining $Z$.
Note that we allow $\Lambda = \emptyset$,
in which case we have $Z = \emptyset$ and $I=(1)$.

We fix $d \geq \max\{2,n\}$.
Then we set
\begin{align}
    &P_i = \sum_{\substack{\alpha = (\alpha_i)_i \in \ZZ_{\geq0}^{n+1} \\ |\alpha|=d, \alpha_i=0 \\ x^\alpha \in I}} x^\alpha,
\end{align}
where we use the standard multi index notation.
We define $g_i \colon \PP^n_\CC \longrightarrow \PP^n_\CC$
by
\begin{align}
    g_i(x_0:\cdots :x_n) = 
    (x_0^d : \cdots : x_{i-1}^d : 
    x_i^d + P_i : x_{i+1}^d : \cdots : x_n^d )
\end{align}
and define 
\begin{align}
    f = g_n \circ \cdots \circ g_0.
\end{align}
Note that $g_i$'s are morphism by the choice of $P_i$
and $g_i^*\OO(1) \simeq \OO(d)$, $f^*\OO(1) \simeq \OO(d^{n+1})$.

\begin{proposition}\label{prop:example-giventotallyinvariantset}
    In the above notation, $Z$ is the 
    maximal proper totally invariant Zariski closed subset
    of $f$. That is, 
    $f^{-1}(Z) = Z$
    and for any closed subset $E \subseteq \PP^n_\CC$ with 
    $f^{-1}(E) = E$, we have $E \subseteq Z$.
\end{proposition}

\begin{lemma}\label{lemma:Zistotallyinvariant}
    For any $0\leq i \leq n$ and $\lambda \in \Lambda$,
    we have $g_i^{-1}(L_{A_\lambda}) = L_{A_\lambda}$.
    In particular, $f^{-1}(Z) = Z$.
\end{lemma}
\begin{proof}
    If $i \not\in A_\lambda$,
    $g_i^{-1}(L_{A_\lambda}) = (x_j^d=0 \mid j \in A_\lambda) = L_{A_\lambda}$ as sets.
    If $i \in A_\lambda$, then
    $g_i^{-1}(L_{A_\lambda}) = (x_j^d=0 \mid j \in A_\lambda \setminus \{i\}) \cap (x_i^d + P_i=0)$.
    Since $P_i=0$ on $L_{A_\lambda}$ and $x_i$ does not appear in $P_i$, $P_i = 0$ on $(x_j=0 \mid j \in A_\lambda \setminus \{i\})$.
    Thus we proved $g_i^{-1}(L_{A_\lambda}) = L_{A_\lambda}$.
\end{proof}

To prove the maximality of $Z$, we prepare a lemma.
For a subvariety $V \subseteq \PP^n_\CC$, we set
\begin{align}
    s(V) = \# \{j \in \{0,\dots, n\} \mid V \subseteq V_+(x_j)\}.
\end{align}
In particular, for a point $a = (a_0:\dots : a_n) \in \PP^n(\CC)$,
$s(a) = \#\{ j \mid a_j =0\}$.

\begin{lemma}\label{lemma:multiplicity-of-g}
    Let $d \geq 1$.
    Let $0 \leq i \leq n$ and 
    $P$ be a homogeneous polynomial in $x_0,\dots, x_n$ of degree $d$
    such that $x_i$ does not appear.
    Let $g \colon \PP^n_\CC \longrightarrow \PP^n_\CC$
    be $g(x_0:\cdots : x_n) = (x_0^d : \cdots : x_i^d + P : \cdots : x_n^d)$.
    Then for any $a \in \PP^n(\CC)$, we have
    \begin{align}
        e_g(a) := \length (\OO_{\PP^n,a} / g^*\gm_{g(a)} \OO_{\PP^n,a})
        = d^{s(a)}.
    \end{align}
\end{lemma}
\begin{proof}
    By rearranging coordinates, we may assume $i = n$.
    First we consider the case where there is $j \neq n$
    such that $x_j(a) \neq 0$.
    By rearranging coordinate again, we may assume $j = 0$.
    In the affine coordinate $z_l = x_l/x_0, 1 \leq l \leq n$,
    $g(z_1,\dots, z_n) = (z_1^d,\dots, z_{n-1}^d, z_n^d + Q)$
    for some polynomial $Q$ in $z_1,\dots, z_{n-1}$.
    Let $a_j$ be the affine coordinates of $a$.
    Thus the ideal $g^*\gm_{g(a)}$ is generated by
    \begin{align}
        &z_l^d - a_l^d \ \text{for $1 \leq l \leq n-1$ and}\\
        & \f := z_n^d - a_n^d + Q - Q(a).
    \end{align}
    If we set $S = \{l \mid \text{$1 \leq l \leq n-1$ and $a_l=0$}\}$, then
    \begin{align}
        \OO_{\PP^n, a}/(z_l^d - a_l^d \mid 1 \leq l \leq n-1)
        \simeq
        \big(\CC[z_l \mid l \in S \cup \{n\}]/(z_l^d \mid l \in S) \big)_{(z_n - a_n, z_l, l \in S)} =: R.
    \end{align}
    Note that the image of $Q - Q(a)$ in $R$
    is contained in the ideal $(z_l \mid l\in S)$.
    Thus we have
    \begin{align}
        \OO_{\PP^n,a}/g^*\gm_g(a)\OO_{\PP^n,a}
        \simeq 
        \big((\CC[t_1,\dots, t_s]/(t_1^d,\dots, t_s^d))[t]/(t^d - a_n^d + h(t_1,\dots,t_{n-1})) \big)_{(t_1,\dots, t_s, t-a_n)}
    \end{align}
    where $s = \# S$ and for some $h$ with constant term zero.
    The length of this ring is $d^{s}$ if $a_n \neq 0$
    and $d^{s+1}$ if $a_n=0$ as we desired.

    Finally, we consider the case where $x_j(a)=0$ for 
    all $1 \leq j \leq n-1$.
    In this case, $g(a) = a$.
    We have
    \begin{align}
        g^*(x_j/x_n) = \frac{1}{1 + P/x_n^d} \frac{x_j^d}{x_n^d}.
    \end{align}
    Since $a + P/x_n^d$ is unit in $\OO_{\PP^n, a}$, 
    we get $e_g(a) = d^n$ as desired.
\end{proof}

\begin{proof}[Proof of \cref{prop:example-giventotallyinvariantset}]
    We proved total invariance of $Z$ in  \cref{lemma:Zistotallyinvariant}.
    Let $E \subsetneq \PP^n_\CC$ be a proper closed subset
    such that $f^{-1}(E)=E$.
    We need to show $E \subseteq Z$.
    Let us set $E_i = g_{i-1} \circ\cdots \circ g_0 (E)$
    for $i=0,\dots, n+1$.
    Here $E_0=E$ by convention and $E_{n+1} = E$ since $f(E)=E$.
    \begin{align}
        E = E_0 \xrightarrow{g_0} E_1 \xrightarrow{g_1} E_2 \xrightarrow{g_2} \cdots \xrightarrow{g_{n-1}} E_n 
        \xrightarrow{g_n} E_{n+1}=E.
    \end{align}
    Since $f^{-1}(E)=E$, we have $g_i^{-1}(E_{i+1}) = E_i$.
    Moreover, $g_i$ induces bijection between irreducible components 
    of $E_i$ and $E_{i+1}.$

    \begin{claim}\label{claim:irrcomp-of-E-coordlin}
        Irreducible components of $E_0,\dots E_n$ are all coordinate linear subvarieties.
    \end{claim}
    \begin{proof}[Proof of \cref{claim:irrcomp-of-E-coordlin}]
    Let $V \subseteq E_{i}$ be an irreducible component 
    and let $W = g_i(V)$:
    \begin{equation}
        \begin{tikzcd}
            &[-2em]\PP^n_\CC \arrow[r, "g_i"]& \PP^n_\CC \\
            g_i^{-1}(E_{i+1}) \arrow[r, phantom, "="]&E_i \arrow[u, phantom, "\subseteq", sloped] \arrow[r]& E_{i+1} \arrow[u, phantom, "\subseteq", sloped]\\
            g_i^{-1}(W) \arrow[r, phantom, "="]&V \arrow[u, phantom, "\subseteq", sloped] \arrow[r, "g_i|_V"] & W \arrow[u, phantom, "\subseteq", sloped]
        \end{tikzcd}
    \end{equation}
    We are applying \cref{lemma:multiplicity-of-g} to $g_i$.
    Notice that for a general point $a \in V$,
    we have $s(a) = s(V)$.
    Thus for general $a \in V$, we have $e_g(a) = d^{s(V)}$.
    Thus we have 
    \begin{align}
        d^n = \deg g_i = \deg (g_i|_V) d^{s(V)}.
    \end{align}
    On the other hand, by the projection formula, we have
    \begin{align}
        \deg(g_i|_V) \deg W = (c_1(\OO(1))^{\dim W} \cap {g_i}_*[V])
        = (c_1(\OO(d))^{\dim W} \cap [V]) = d^{\dim V} \deg V.
    \end{align}
    Thus we have
    \begin{align}
        \deg W = d^{\dim V - n + s(V)} \deg V 
        = d^{s(V) - \codim_{\PP^n} V}\deg V.
    \end{align}
    Note that $s(V) \leq \codim_{\PP^n} V$,
    and $s(V) \leq \codim_{\PP^n} V$ if and only if
    $V$ is a coordinate linear subvariety.
    If we apply $g_i, g_{i+1},\dots$ repeatedly, then
    the image of $V$ eventually back to $V$ itself.
    Thus we must have $s(V)=\codim_{\PP^n}V$,
    hence $V$ is a coordinate linear subvariety.
    \end{proof}

    \begin{claim}\label{claim:gicoordlinsub}
        Let $\emptyset \neq A \subseteq \{0,\dots, n\}$.
        If $g_i(L_A)$ is a coordinate linear subvariety,
        then $g_i(L_A)=L_A$.
    \end{claim}
    \begin{proof}[Proof of \cref{claim:gicoordlinsub}]
        If $i \not\in A$, then we have $g_i(L_A) \subseteq L_A$
        and hence $g_i(L_A) = L_A$.
        Suppose $i \in A$.
        If $P_i$ is zero on $L_A$, then $g_i(L_A) \subseteq L_A$
        and again we get $g_i(L_A) = L_A$.
        Suppose $P_i$ is not constantly zero on $L_A$.
        Then $g_i(L_A) \subseteq L_{A \setminus \{i\}}$,
        but for any $j \not\in A \setminus \{i\}$,
        $x_j$ does not vanish on $g_i(L_A)$.
        Thus $g_i(L_A)$ is not a coordinate linear subvariety in this case.
    \end{proof}

    Let us pick an arbitrary irreducible component of $E$.
    By \cref{claim:irrcomp-of-E-coordlin}, we can write it 
    as $L_A$ for some non-empty $A \subseteq \{0,\dots, n\}$.
    By \cref{claim:irrcomp-of-E-coordlin} and \cref{claim:gicoordlinsub}, $g_i(L_A) = L_A$ for all $i$.
    In particular, for $i \in A$, $P_i = 0$ on $L_A$ constantly.
    
    Suppose $L_A \not\subseteq Z$.
    Let $B = \{0 ,\dots, n\} \setminus A$.
    Since $L_A \not\subseteq L_{A_\lambda}$ for all
    $\lambda \in \Lambda$, we have $A_\lambda \cap B \neq \emptyset$.
    Hence the monomial $\prod_{j \in B}x_j$ vanishes on $Z$.
    Pick a $b \in B$ and set
    \begin{align}
        M = x_b^{d-\#B}\prod_{j \in B} x_j.
    \end{align}
    Note that the exponent is non-negative since $\# B \leq n \leq d$.
    Thus $M \in I$. 
    For any $i \in A$, $M$ is a summand of $P_i$
    since $x_i$ does not appear in $M$.
    As $M|_{L_A}$ is not constantly zero, 
    we have $P_i$ is not constantly zero on $L_A$, 
    a contradiction.
    (Note that two monomials in $P_i$ do not coincide 
    after restricted on $L_A$ unless both of them vanish on $L_A$.)
\end{proof}

\appendix

\section{Bott vanishing on projective spaces}

In this section, we work over a field $k$.
Let $\PP^N_k$ be the projective space with the structure sheaf $\OO$.
Set $W = H^0(\PP^N_k, \OO(1))$, $\Omega = \Omega_{\PP^N_k/k}$.
Recall the Euler sequence 
\begin{align}
    0 \longrightarrow \Omega \longrightarrow \mathcal{F} \xlongrightarrow{\epsilon} \OO \longrightarrow 0
    \label{eulerseq}
\end{align}
where $\mathcal{F} = W \otimes_k \OO(-1)$.
For $r \geq 1$, the $\OO$-morphism $\epsilon$ induces 
\begin{align}
    \delta_r \colon \bigwedge^r \mathcal{F} \longrightarrow \bigwedge^{r-1} \mathcal{F}
\end{align}
by the following rule:
\begin{align}
    \delta_r(\alpha_1 \wedge \cdots \wedge \alpha_r)
    = 
    \sum_{l=1}^r (-1)^{l-1} \epsilon(\alpha_l) \alpha_1 \wedge \cdots
    \wedge \widehat{\alpha_l} \wedge \cdots \wedge \alpha_r.
\end{align}
Since \cref{eulerseq} splits locally, we can see that 
the image of $\delta_r$ is $\Omega^{r-1}$ and the kernel is
$\Omega^r$:
\begin{align}
    0 \longrightarrow \Omega^r \longrightarrow \bigwedge^r \mathcal{F}
    \xlongrightarrow{\delta_r}  \Omega^{r-1} \longrightarrow 0
    \label{wedgedeulerseq}
\end{align}
is exact.

\begin{proposition}\label{appendix-prop:bottvanforprojsp}
    Let $r \geq 1$.
    For $e \in \ZZ$, we write $\delta_r(e)$ the $e$-th twist of $\delta_r$.
    Then the induced morphism
    \begin{align}
        H^0(\PP^N_k, \bigwedge^{r+1} \mathcal{F}(e)) \longrightarrow
        H^0(\PP^N_k, \Omega^{r}(e))
    \end{align}
    is surjective.
    Moreover, we have
    \begin{align}
        H^0(\PP^N_k, \Omega^d(d)) = 0
    \end{align}
    for all $d \geq 1$, and
    \begin{align}
        H^{q}(\PP^N_k, \Omega^d(m)) = 0
    \end{align}
    for all $q, m \geq 1$ and $0 \leq d \leq N$.
\end{proposition}
\begin{proof}
    We set $S = \Sym W$. Its degree $l$ part is written by $S_l$.
    By the preceding argument, we have
    \begin{equation}
        \begin{tikzcd}
            &\Omega^r(e) \arrow[d, hook]&\\
            \bigwedge^{r+1}\mathcal{F}(e) \arrow[r,"\delta_{r+1}(e)"] \arrow[ru]& \bigwedge^{r}\mathcal{F}(e)
            \arrow[r, "\delta_{r}(e)"] &\bigwedge^{r-1}\mathcal{F}(e)
        \end{tikzcd}
    \end{equation}
    where $\ker \delta_r(e) = \Omega^r(e)$.
    By taking global sections, we get
    \begin{equation}
        \begin{tikzcd}
            &H^0(\PP^N_k, \Omega^r(e)) \arrow[d, hook]&\\
            \bigwedge^{r+1}W \otimes_k S_{e-r-1} \arrow[r,"\delta_{r+1}(e)"] \arrow[ru]& \bigwedge^{r}W \otimes_k S_{e-r}
            \arrow[r, "\delta_{r}(e)"] &\bigwedge^{r-1}W \otimes_k S_{e -r +1}
        \end{tikzcd}.\label{keydaigramobtainedbyeulerseq}
    \end{equation}
    We claim that the bottom row is exact.
    Indeed, if we define
    \begin{equation}
        \begin{tikzcd}
            \bigwedge^r W \otimes_k S  \arrow[r, "\partial" ]&
            \bigwedge^{r-1} W \otimes_k S\\[-1ex]
            (w_1 \wedge \cdots \wedge w_r) \otimes f \arrow[r,mapsto]& 
            \sum_{l=1}^r(-1)^{l-1} w_1 \wedge \cdots \wedge
            \widehat{w_l} \wedge \cdots \wedge w_r \otimes w_l f
        \end{tikzcd}
    \end{equation}
    the sequence 
    \begin{align}
        0 \to \bigwedge^{N+1} W \otimes_k S \xrightarrow{\partial}
        \bigwedge^{N}W \otimes_k S  \xrightarrow{\partial} 
        \cdots \xrightarrow{\partial} \bigwedge^2W \otimes_k S
        \xrightarrow{\partial} W \otimes S
        \xrightarrow{\partial} S \to k
        \label{koszcpx}
    \end{align}
    is the Koszul complex of $S$ and a $k$-basis of $W$,
    hence exact.
    If we set grading on $\bigwedge^r W \otimes_k S$
    by $\bigwedge^r W \otimes_k S 
    = \bigoplus_{e \in \ZZ} \bigwedge^rW \otimes_k S_{e -r}$,
    the maps $\partial$ are all degree preserving.
    The degree $e$ part of $r-1, r, r+1$-th terms of 
    \cref{koszcpx} is exactly equal to the bottom row of
    \cref{keydaigramobtainedbyeulerseq}. This proves our claim.

    The first assertion directly follows from this exactness
    and \cref{keydaigramobtainedbyeulerseq}.
    Moreover, if we set $r = d , e = d$,
    then $H^0(\bigwedge^{d+1}\mathcal{F}(d)) = 0$.
    This implies $H^0(\Omega^d(d)) = 0$ for $d \geq 1$.
    
    Finally, by (\ref{wedgedeulerseq}) and our first assertion, 
    we have the following 
    long exact sequence for $r \geq 1$, $m \geq 1$:
    \begin{equation}
    \begin{aligned}
        0 &\longrightarrow H^1(\Omega^r(m))
        \longrightarrow H^1({\textstyle\bigwedge^r}\mathcal{F}(m))
        \longrightarrow H^1(\Omega^{r-1}(m)) \\
        &\longrightarrow H^2(\Omega^r(m))
        \longrightarrow \cdots \\
        &\longrightarrow H^N(\Omega^r(m))
        \longrightarrow H^N({\textstyle\bigwedge^r}\mathcal{F}(m))
        \longrightarrow H^N(\Omega^{r-1}(m))
        \longrightarrow 0.
    \end{aligned}
\end{equation}
    
    If $1 \leq r \leq N+1$, all the higher cohomologies of $\bigwedge^r \mathcal{F}(m)$ are zero.
    Starting form $r=N+1$ to $r = 1$, we can inductively see
    that $H^q(\Omega^d(m))=0$ for $q \geq 1, m\geq 1$ and $0 \leq d \leq N$.
\end{proof}

\section{Degree one subvariety is linear}

\begin{lemma}\label{lem:degoneislin}
    Let $Z \subseteq \PP^n_\CC$ be an irreducible reduced closed
    subscheme of degree one. Then $Z$ is a linear subvariety.
\end{lemma}
\begin{proof}
    We may assume $\dim Z =: d \geq 2$.
    We claim that for any line $l \subseteq \PP^n_\CC$,
    we have
    \begin{align}
        l \cap Z = 
        \begin{cases}
            \emptyset\\
            \text{one point as a scheme}\\
            l
        \end{cases}
    \end{align}
    To prove this, assume $l \cap Z \neq l$.
    Then we can find hyperplanes $H_1,\dots, H_d \subseteq \PP^n_\CC$ with the following properties inductively:
    for $i=0,\dots, d-1$,
    $Z \cap H_1 \cap \cdots \cap H_i$ is pure dimension $d-i$
    and $H_{i+1}$ restricts on ${Z \cap H_1 \cap \cdots \cap H_i}$ as an effective Cartier divisor;
    $l \cap Z \subseteq Z \cap H_1 \cap \cdots \cap H_d$ as 
    closed subscheme of $\PP^n_\CC$.
    Then by \cite[Lemma 1.7.2]{Ful98}, we have
    \begin{align}
        1 = (c_1(\OO(1))^d \cap [Z]) =  \deg[Z \cap H_1 \cap \cdots H_d].
    \end{align}
    Thus $l \cap Z$ must be $\emptyset$ or a single point
    as a scheme.

    By this claim, any line passing through two distinct points
    of $Z$ are contained in $Z$.
    Thus the affine cone of $Z$ is a linear subspace.
\end{proof}

\begingroup
\raggedright
\bibliographystyle{amsplain}
\bibliography{totally_invariant}
\endgroup

\end{document}